\documentclass{gtmon_a}
\pdfoutput=1
\usepackage{pinlabel}
\proceedingstitle{Groups, homotopy and configuration spaces (Tokyo
  2005)}
\conferencestart{5 July 2005}
\conferenceend{11 July 2005}
\conferencename{Groups, homotopy and configuration spaces,
                in honour of Fred Cohen's 60th birthday}
\conferencelocation{University of Tokyo, Japan}

\editor{Norio Iwase}
\givenname{Norio}
\surname{Iwase}

\editor{Toshitake Kohno}
\givenname{Toshitake}
\surname{Kohno}

\editor{Ran Levi}
\givenname{Ran}
\surname{Levi}

\editor{Dai Tamaki}
\givenname{Dai}
\surname{Tamaki}

\editor{Jie Wu}
\givenname{Jie}
\surname{Wu}

\title{Automorphic functions for a Kleinian group}                    

\author{Masaaki Yoshida}
\givenname{Masaaki}
\surname{Yoshida}
\address{Department of Mathematics\\Kyushu University\\\newline
Fukuoka 810-8560\\Japan}
\email{myoshida@math.kyushu-u.ac.jp}
\urladdr{}

\dedicatory{Dedicated to Professor Fred
Cohen on his sixtieth birthday}

\volumenumber{13}
\issuenumber{}
\publicationyear{2008}
\papernumber{22}
\startpage{483}
\endpage{497}

\doi{}
\MR{}
\Zbl{}

\arxivreference{}


\received{24 November 2005}
\revised{12 December 2006}
\accepted{19 January  2007}
\published{19 March 2008}
\publishedonline{19 March 2008}
\proposed{}
\seconded{}
\corresponding{}
\version{}

\let\xysavmatrix\xymatrix
\def\xymatrix{\disablesubscriptcorrection\xysavmatrix}
\AtBeginDocument{\let\bar\wbar\let\tilde\wtilde\let\hat\what}

\makeop{home}
\newcommand{\we}{\smash{\rlap{\kern 6pt
\raise 4pt\hbox{\footnotesize $\sim$}}}\longrightarrow}
\makeautorefname{subsection}{subsection}
\makeautorefname{subsubsection}{subsection}
\newcommand{\B}{{\mathbb B}}

\renewcommand{\P}{{\mathbb P}}
\renewcommand{\H}{{\mathbb H}}
\newcommand{\RR}{{\mathbb{R}}}
\def\overset#1#2{{\buildrel #1\over#2}}

\newcommand{\G}{\Gamma}

\def\ip{{1{+}i}}

\begin{document}

\begin{abstract}  
In the paper `Automorphic functions for a Whitehead-complement group'
 \cite{MNY}, Matsumoto, Nishi and Yoshida constructed automorphic
 functions on real 3--dimen\-sional hyperbolic space for a Kleinian
 group called the Whitehead-link-complement group.  For a Kleinian
 group (of the first kind), no automorphic function/form has been
 studied before. In this note, their motivation is presented with a
 historical background.
\end{abstract}

\begin{webabstract}  
In the paper `Automorphic functions for a Whitehead-complement group',
 [Osaka J Math 43 (2006) 63--77] Matsumoto, Nishi and Yoshida
 constructed automorphic functions on real 3--dimensional
 hyperbolic space for a Kleinian group called the
 Whitehead-link-complement group. For a Kleinian group (of the first
 kind), no automorphic function/form has been studied before. In this
 note, their motivation is presented with a historical background.
\end{webabstract}

\begin{asciiabstract}  
In the paper `Automorphic functions for a Whitehead-complement group',
 [Osaka J Math 43 (2006) 63-77] Matsumoto, Nishi and Yoshida
 constructed automorphic functions on real 3-dimensional
 hyperbolic space for a Kleinian group called the
 Whitehead-link-complement group. For a Kleinian group (of the first
 kind), no automorphic function/form has been studied before. In this
 note, their motivation is presented with a historical background.
\end{asciiabstract}

\maketitle

\section{Introduction}
In \cite{MNY} Matsumoto, Nishi and Yoshida constructed automorphic
functions on real 3--dimensional hyperbolic space for a Kleinian group
called the Whitehead-link-complement group.  For a Kleinian group (of
the first kind), no automorphic function/form has been studied before,
except by Matsumoto and Yoshida \cite{MY}.

In this note, we first recall (from a very elementary level) branched
covers of the complex projective line, especially those with three
branch points; these are prototypes of our branched cover story.  We
also define, as higher dimensional generalizations, several branched
covers of complex projective spaces, when the universal ones are the
complex ball (complex hyperbolic sapce) and we briefly mention the
history.

The main topic of this note is another higher dimensional
generalization: branched covers of the 3--sphere.  One of the simplest
links/knots whose complements admit hyperbolic structures is the
Whitehead link $L$; there is a discrete subgroup $W$ of the
automorphism group of real hyperbolic 3--space $\H^3_\RR$, such
that the quotient space $\H^3_\RR/W$ is homeomorphic to the complement
$S^3-L$. We construct automorphic fucntions on $\H^3_\RR$ for $W$ and
make the projection $\H^3_\RR\to S^3-L$ explicit.

\section{Covers of the complex projective line} 
Let $X=\P^1$ be the complex projective line, which is also called the
Riemann sphere. We are interested in its coverings. Since $X$ is
simply connected, there is no non-trivial covering unless we admit
branch points.

\subsection{One branch point}
If we admit only one branch point, which we can assume to be the point
at infinity, then since $X-\{\infty\}\cong\C$ is still simply
connected, there is no non-trivial cover.

\subsection{Two branch points}
If we admit two branch points, which we can assume to be the origin
and the point at infinity, and assign $p,q\in\{2,3,\cdots,\infty\}$ as
the indices, respectively, then such a covering of $X$ exists if and
only if $p=q$. Let us denote such a {\it covering space} by $Z$, and
the {\it projection} by $\pi\co Z\to X$; its (multi-valaued) inverse
is called the {\it developing map} $s\co X\to Z$. The covering (deck)
transformation group, which is a group of automorphisms of $Z$, is
denoted by $\Gamma$; note that the projection $\pi$ is
$\Gamma$--automorphic, ie, invariant under the action of $\Gamma$.
\subsubsection{Case $p=q<\infty$}
The covering space and the projection are given by
$$\pi\co \P^1\cong Z\ni z\longmapsto x=z^p\in X,$$
and the developing map is the multi-valed map $x\mapsto z=x^{1/p}$.
Note that the projection is invariant under the finite group
$$\Gamma=\{z\mapsto e^{\frac{2\pi ik}p}z\mid k=1,\dots,p\}.$$ 
\subsubsection{Case $p=q=\infty$} 
The covering space and the projection are given by
$$\pi\co \C\cong Z\ni z\longmapsto x=\exp z\in X$$ and the developing
map is the multi-valed map $x\mapsto z=\log x$.  Note that the
projection is invariant under the infinite group
$$\Gamma=\{z\mapsto z+2\pi ik\mid k\in \Z\}.$$ 
In this way, two points on a sphere --  a simple geometric object -- 
naturally leads to very important functions the exponential and the logarithm functions.
\subsection{Three branch points}
If we admit three branch points, which we can assume to be $x=0,1,\infty$,
and assign $p,q, r\in\{2,3,\cdots,\infty\}$ as the indices,
respectively, then such a covering of $X$ always exists; there are
many such. Let $Z$ be the biggest one, the universal branched covering
with the pre-assigned indices; this is characterized as the simply
connected one. There are only three simply connected 1--dimensional
complex manifolds. The ramification indices determine the nature of
$Z$. We tabulate the three cases with familiar names for the deck
transformation groups.
$$\begin{array}{ccccccc}
\mbox{type}&&\frac1p+\frac1q+\frac1r&&    Z  &&\mbox{deck transformation group }\Gamma\\[3mm]
\mbox{elliptic}&&>1                 && \P^1  &&\mbox{polyhedral groups}\\[3mm]
\mbox{parabolic}&&=1   &&\C&&\mbox{1-$\dim_\C$ crystallographic groups}\\[3mm]
\mbox{hyperbolic}&&<1&&\H&&\mbox{triangle (Fuchsian) groups}\end{array}$$  
Complex crystallographic groups of dimension $n$ are by definition, groups of affine transformations of $\C^n$ with compact quotients. A Fuchsian group (of the first kind) is by definition a discrete subgroup of $SL(2,\R)$ acting on the upper half-space
$$\H=\{\tau\in\C\mid\Im\tau>0\},$$
such that a(ny) fundamental domain has finite volume.
\subsubsection{Developing maps}
In each case the developing map is given by the Schwarz map
$$s\co X\ni x\mapsto z=u_0(x)/u_1(x)\in Z,$$
defined by the ratio of two (linealy independent) solutions of the hypergeometric differential equation $E(a,b,c)$:
$$x(1-x)u''+\{c-(a+b+1)x\}u'-abu=0,$$
where the parameters $a,b,c$ are determined by the condition
$$|1-c|=\frac1p,\quad|c-a-b|=\frac1q,\quad|a-b|=\frac1r.$$  See Yoshida
\cite{Yos} for hypergeometric functions, thir fundamental properties
and the Schwarz map.
\subsubsection{Projections}\label{sec:2.3.2}
In each elliptic case, the projection is given by an invariant of a
polyhedral group $\Gamma$, which is classically known. 

In each parabolic case, the projection is given by an elliptic function.

However in hyperbolic cases, only for a few examples, are the
projections explicitly constructed.  Here we explain a typical example
when $(p,q,r)=(\infty,\infty,\infty)$. The developing map is the
Schwarz map of the equation $E(1/2,1/2,1)$, the deck transformation
group is (conjugate to) the congruence subgroup
$$\Gamma(2)=\{g\in SL(2,\Z)\mid g\equiv \mbox{identity mod 2}\}$$
of the elliptic full modular group $SL(2,\Z)$.
The projection is given by the lambda function
$$\lambda(\tau)=\left(\frac{\vartheta_{01}(\tau)}{\vartheta_{00}(\tau)}\right)^4,\quad \tau\in\H,$$
where
$$\vartheta_{00}(\tau)=\sum_{n\in\Z}q^{2n^2},\quad  
\vartheta_{01}(\tau)=\sum_{n\in\Z}(-1)^nq^{2n^2},\qquad q=\exp\pi i \tau/2.$$
\subsubsection{Note} 
Three points on a sphere --  a simple geometric object -- naturally leads to 
interesting kind of mathematics such as polyhedral groups and their invariants, Fuchsian automorphic forms/functions and the hypergeometric differential equation/functions.
\subsection{Four or more branch points} 
If we admit four branch points, which we can assume $x=0,1,\infty,t$,
and assign any four indices ${\bf p}=(p_0,\cdots,p_3)$, then such a
covering of $X$ always exists. Let $Z$ be the biggest one. If the
indices are ${\bf p}=(2,\cdots,2)$ then $Z\cong\C$, otherwise $Z\cong
\H$. The developing map is the Schwarz map of a second-order Fuchsian
equation $E(t,{\bf p})$. The $t$--dependence of the coefficients of
this equation is difficult to analyse, and when $t$ is in a generic
position, the coefficients have no explicit expressions. The
projection is known only for sporadic values of $t$ and ${\bf p}$.
Though there have been many attempts for tackling this difficulty the
goal is still far away.  I am afraid that `four points on $\P^1$' is
an object too difficult for human beings.  There is no hope for more
than four points.
\section{Covers of the complex projective spaces}
When one encounters a serious difficulty, one of the ways to proceed is to
turn to high dimensional analogues. 
\subsection{Covers of $\P^2$}
The first attempt was made by
E~Picard \cite{Pic}. He considers the six lines
$$x_0x_1x_2(x_0-x_1)(x_1-x_2)(x_2-x_0)=0$$
on the projective plane $X(\cong \P^2)$ with homogeneous coordinates
$x_0:x_1:x_2$, and studies the universal branched covering $Z$ of $X$
ramifying along the six lines with indices 3. It turns out that $Z$ is
isomorphic to the complex 2--dimentional ball
$$\B^2=\{z_0:z_1:z_2\in\P^2\mid |z_0|^2-|z_1|^2-|z_2|^2>0\},$$
and that the developing map is the Schwarz map 
$$s\co X\ni x\longmapsto z=u_0(x):u_1(x):u_2(x)\in Z,$$
defined by linearly independent solutions $u_0,u_1$ and $u_2$ of 
Applell's hypergeometric differential equation $F_1$ (with special
parameters), a system of linear partial differential equations defined 
on $X$ of rank three with singularities along the six lines. The
deck transformation group is an arithmetic subgoup of the group of
automorphisms of $\B^2$ defined over $\Z(1^{1/3})$, 
and the projection can be expicitly expressed in terms of Riemann theta
functions (this part is completed by H~Shiga \cite{Shi}).
\subsection{Covers of $\P^n$}
T~Terada \cite{Ter1,Ter2} considers the hyperplanes
$$\prod_{i=1}^nx_i\cdot\prod_{i,j=1,i\not=j}^n(x_i-x_j)=0$$
on the projective $n$--space $X(\cong \P^n)$ with homogeneous coordinates
$x_0:\cdots:x_n$, and studies the universal branched covering $Z$ of $X$
ramifying along these hyperplanes with different indices. It turns out, 
for some choices (finite possibilities) of indices, that $Z$ is
isomorphic to the complex $n$--dimentional ball
$$\B^n=\{z_0:\cdots:z_n\in\P^n\mid |z_0|^2-\sum_{i=1}^n|z_n|^2>0\},$$
and that the developing map is the Schwarz map 
$$s\co X\ni x\longmapsto z=u_0(x):\cdots:u_n(x)\in Z,$$
defined by linearly independent solutions $u_0,\cdots, u_n$ of 
Applell's hypergeometric differential equation $F_D$ (with 
parameters corresponding to the indices chosen), 
a system of linear partial differential equations defined 
on $X$ of rank $n+1$ with singularities along those hyperplanes. The
deck transformation group is a discrete subgoup (not necessarily
arithmetic, see Deligne and Mostow \cite{DM}) of the group of
automorphisms of $\B^n$. 
For several cases, the projection can be explicitly expressed in terms 
of Riemann theta functions (Shiga, \cite{Shi},  Matsumoto and Terasoma
\cite{MT1,MT2}).
\subsection{Note}
Recently further studies have been made; the objects and techniques
require considerably more algebraic geometry and representation
theory.  For more detail see Allcock, Carlson and Toledo \cite{ACT},
or Couwenberg, Heckman and Looijenga \cite{CHL}.
\section{Covers of the 3--sphere}
If one recalls that the complex projective line is also called the 
(Riemann) sphere and that the Poincar\'e upper half-plane $\H$ is 
just the real hyperbolic 2--space, 
another high-dimensional generalization would be made by
the 3--dimensional sphere $S^3$ and the real hyperbolic 3--space
$$\H^3_\RR=\{\mbox{$2\times2$ positive Hermitian matrices}\}/{\R_{>0}}.$$
Since the branch locus must be a submanifold of codimension 2, 
we consider links/knots as branch loci.  

If our branch locus is a trivial knot, nothing interesting happens.

If our branch locus is the trefoil knot, several interesting things 
happen, but the universal cover (branched cover with branch index 
$=\infty$) is not $\H^3_\RR$.

If our branch locus is the Whitehead link or the figure eight knot,
it is known that its universal cover is $\H^3_\RR$. I believe that
these are the simplest ones, and that these are the ones we 
should/could study. For more complicated links/knots, their complements would permit hyperbolic structures, but one can not expect nice mathematics from 
the view point of function theory. In this note, we just treat the
Whitehead link.
\subsection{Whitehead-link-complement group}
It is known that the complement of the Whitehead link $L$
admits a hyperbolic structure; 
that is, the universal cover of $S^3-L$ is $\H^3_\RR$, 
or there is a discrete subgroup 
$W$ of the automorphisms of $\H^3_\RR$,
such that $S^3-L$ is {\it homeomorphic} 
to the quotient space $\H^3_\RR/W$.

The group of automorphisms of $\H^3_\RR$  is generated by $PGL(2,\C)$
acting as
$$GL(2,\C)\ni g\co \H^3_\RR\ni z\longmapsto gz^t\bar g\in\H^3_\RR,$$
and the (orientation reversing) transpose operation
$$T\co \H^3_\RR\ni z\longmapsto {}^tz\in\H^3_\RR.$$
The Whitehead-link-complement group $W$ is generated by the two
elemnts
$$\left(\begin{array}{cc}1&i\\0&1
\end{array}\right)\quad\mbox{and}\quad
\left(\begin{array}{cc}1&0\\\ip &1
\end{array}\right).$$ 
This group is a subgroup of the principal congruence subgroup $\Gamma(1+i)$ of the full modular group 
$$\Gamma=GL(2,\Z[i])$$ of
finite index. It is not normal; the conjugate one is $\bar W$ (complex
conjugate).  We have the homeomorphism
$$\home\co \H^3_\RR/W\ \sim\ S^3-L.$$
\subsection{What should be studied?}
\begin{figure}[ht!]
\begin{center}
\includegraphics[scale=.9]{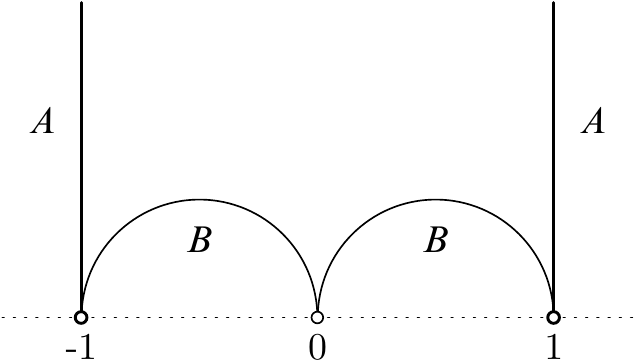}% 
\end{center}
\vspace*{-.5cm}\caption{A fundamental domain for $\Gamma(2)$\label{Gamma2}}
\end{figure}

We stated in \fullref{sec:2.3.2} that we have an isomorphism (of complex analytic 
varieties) $$\lambda\co \H/\Gamma(2)\quad\overset\cong\longrightarrow
\quad\P^1-\{0,1,\infty\}.$$ 
If it is just a homeomorphism
$$\H^2_\RR/\Gamma(2)\sim S^2-\{0,1,\infty\},$$ the picture is easy to
see: a schoolchild could glue the corresponding sides of the
fundamental domain of $\Gamma(2)$ (see \fullref{Gamma2}) and obtain a
balloon with three holes.  The explicit expression (in terms of the
theta functions) of $\lambda$ makes this isomorphism into interesting
mathematics.

\begin{figure}[ht!]
\cl{\includegraphics[width=7cm]{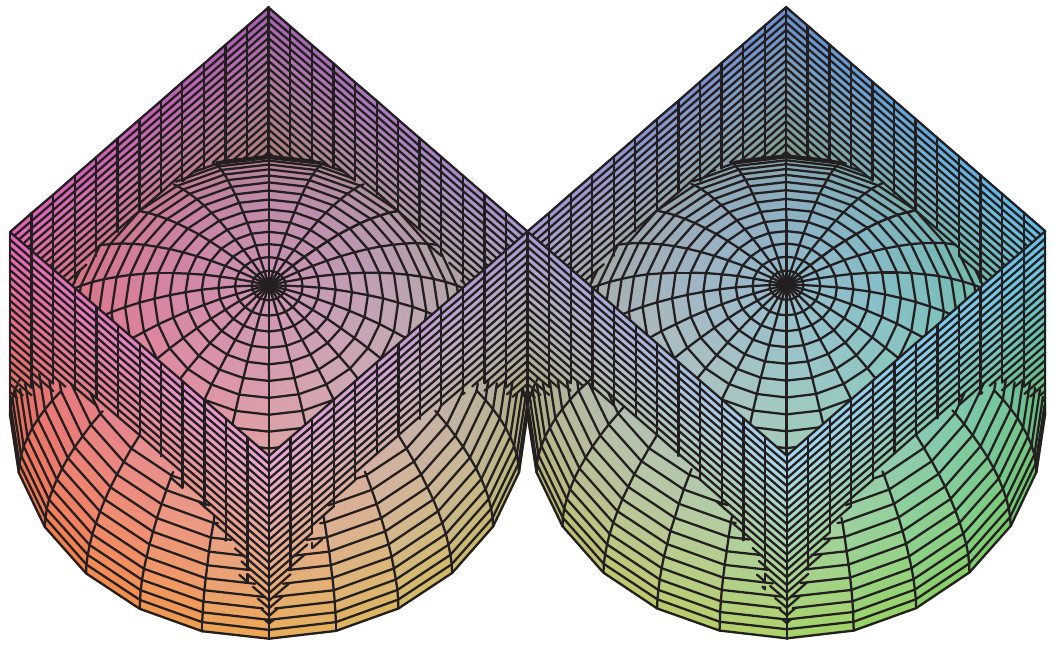}\quad
\includegraphics[scale=.8]{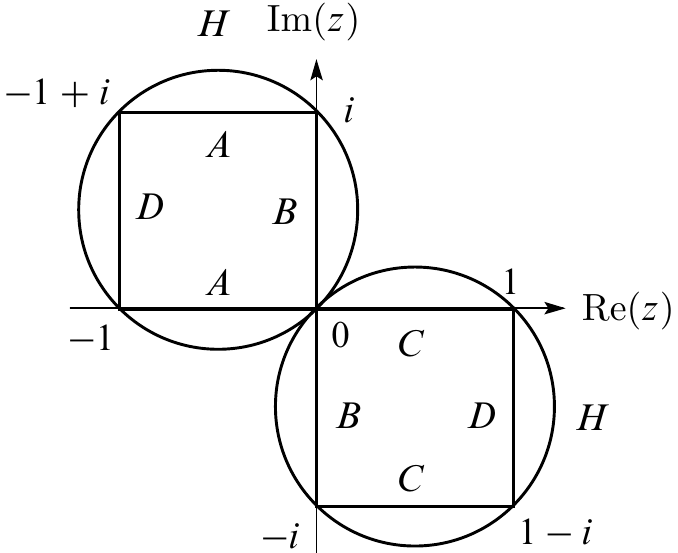}}
\caption{Fundamental domain $FD$ of $W$ in $\H^3_\RR$}
\label{fundW}
\end{figure}

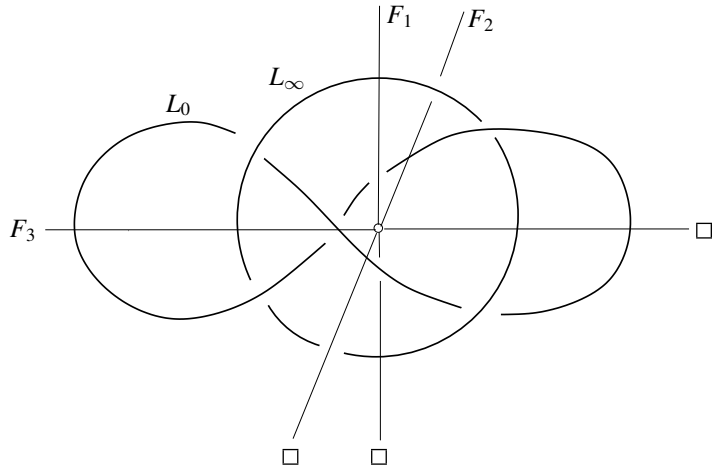
\begin{figure}[ht!]
\begin{center}
\labellist\small
\pinlabel $F_1$ [l] at 359 456
\pinlabel $F_2$ [l] at 446 452
\pinlabel $F_3$ [r] at 0 226
\pinlabel $L_\infty$ [br] at 287 372
\pinlabel $L_0$ [b] at 144 341
\pinlabel $\sq$ [t] at 264 0
\pinlabel $\sq$ [t] at 359 0
\pinlabel $\sq$ [l] at 691 226
\endlabellist 
\includegraphics[width=.65\hsize]{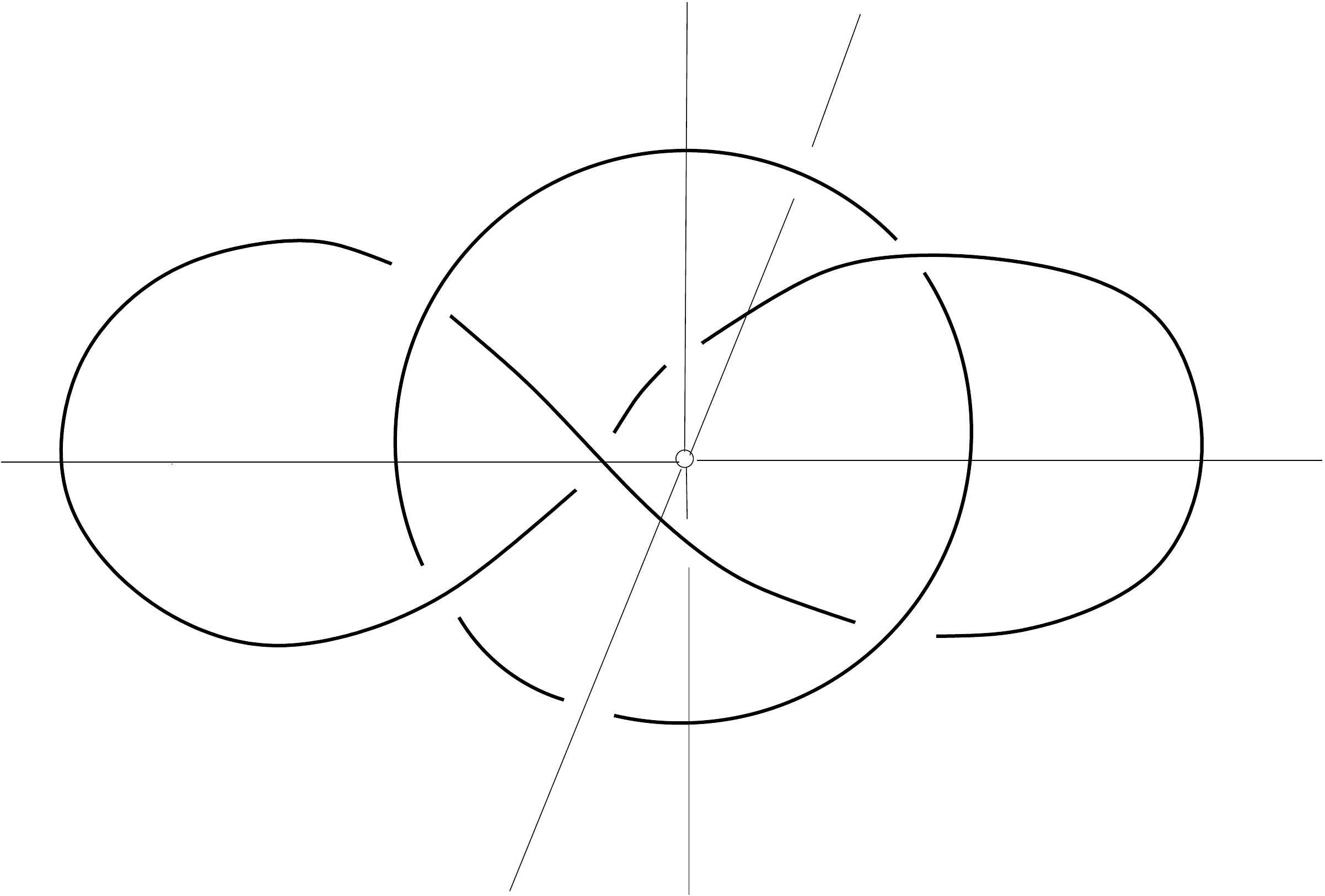}
\end{center}
\caption{Whitehead link with its symmetry axes\label{figwhlink}}
\end{figure}

A fundamenatal domain (consisting of two pyramids) of the
Whitehead-link-comple\-ment group $W$ is shown in \fullref{fundW} (cf
Wielenberg \cite{Wie}); here hyperbolic 3--space is realized as
the upper half-space model
$$\{(z,t)\in\C\times\R\mid t>0\}.$$
The faces of the two pyramids are patched as follows: the eight walls with the same labels ($A,B,C,D$) are patched together, and the two hemi-spheres labelled $H$ are also patched.
The group $W$ has two cusps. They are represented by the vertices of the pyramids:
$$(z,t)=(*,+\infty),\quad (0,0)\sim(\pm i,0)\sim(\pm1,0)\sim(\mp 1\pm i,0).$$
Though it is not so easy as in the 2--dimensional case above to see that the 
fundamental domain modulo this patching is homeomorphic to the complement
of the  Whitehead link $L=L_0\cup L_\infty$ (in
$S^3=\R^3\cup\{\sq\}$)
shown in \fullref{figwhlink},  it is still not advanced mathematics.

Our goal is to find functions on $\H^3_\RR$ invariant under the action of 
$W$, and make the projection
$$\pi\co \quad \H^3_\RR\ \longrightarrow\ S^3-L$$ 
(which induces $\home\co \H^3_\RR/W\ \we  S^3-L$) explicit. 

Though the automorphic forms/functions on the Poincar\'e upper
half-plane $\H$ and those on their high complex-dimensional generalizations
(such as Siegel upper half-spaces, complex balls, and other Hermitian
symmetric spaces) have been studied in detail, those on $\H^3_\RR$
have never been studied before; this is really surprising.

As `function theory' usually means `theory of complex analytic
functions', most (explicit) functions are defined naturally in complex
domains. So we had to rely first on a Hermitian (complex anlytic)
setting, and then to restrict to a real subvariety isomorphic to
$\H^3_\RR$.
\subsection{A Hermitian setting}
More than fifteen years ago, we made in \cite{MSY} a generalization of the 
story (explained above) of the isomorphism
$$\lambda\co  \H/\Gamma(2)
\quad\overset\cong\longrightarrow\quad\P^1-\{0,1,\infty\}.$$
The developing map $\P^1-\{0,1,\infty\}\to\H^1$ can be considered as 
the period map of the family
of elliptic curves
$$t^2=(s-x_1)\cdots(s-x_4),$$
where $(x_1,\dots,x_4)$ is a quadruple of distict points in $\P^1$,
which can be, modulo $PGL(2,\C)$ action, normalized as
$$(0,1,\infty,x), \quad x\in \P^1-\{0,1,\infty\}.$$
We consider the family of K3 surfaces
$$t^2=\ell_1(s_1,s_2)\cdots\ell_6(s_1,s_2),$$
where $\ell_j\ (j=1,\dots,6)$ are linear forms in 
$(s_1,s_2)\in\P^2$ in general position; modulo $PGL(3,\C)$ action,
one can check that this family is 4--dimensional. Let us call this
parameter space $X$, that is,
$$X=PGL(3,\C)\backslash\left\{(\ell_1,\cdots,\ell_6)\mid \mbox{six ordered lines 
in $\P^2$ in general position}\right\}.$$
The periods of these surfaces -- these are also solutions of the generalized hypergeometric differential equation of type $(3,6)$ defined in $X$ --  defines a map from $X$ to the $2\times2$
upper half-space
$$\H_{2\times2}=\left\{\tau:\ 2\times2\mbox{ complex matrix }\mid
\frac{\tau-{}^t\bar\tau}{2i}>0\right\};$$
the deck transormation group (the monodromy group) is the principal
congruence subgroup 
$$\Delta(1+i)=\{g\in\Delta\mid g\equiv I_4\ \mod\ (1+i)\}$$ 
of the full modular group
$$\Delta=GL(4,\Z[i])\cap\{g\in GL(4,\C)\mid gJ_4\ {}^t\bar g=J_4\},
\quad\mbox{where}\quad J_4=\left(\begin{array}{cc}0&-I_2\\ I_2&0\end{array}\right).$$
The group $\Delta$ acts on $\H_{2\times2}$ as
$$g\cdot\tau=(A\tau+B)(C\tau+D), \quad
g=\left(\begin{array}{cc}A&B\\C&D\end{array}\right)\in\Delta, \ \tau\in\H_{2\times2}.$$
The projection $\pi\co \H_{2\times2}\rightarrow X $ can be expressed in terms of the theta functions
$$\theta{a\choose b}(\tau)=\sum_{n\in\Z[i]^2}\exp\pi 
i\{(n+a)\tau\ {}^t\overline{(n+b)}+2\Re(n\cdot {}^t\bar b)\}$$
with characteristics
$$a,b\in\left(\frac{\Z[i]}{1+i}\right)^2.$$
The projection $\pi$ induces the isomorphism
$\H_{2\times2}/\Delta(1+i)\cong X.$
\subsection{Restriction to $\H^3_\RR$ in $\H_{2\times2}$}
The real hyperbolic 3--space $\H^3_\RR$ lives in $\H_{2\times2}$ as
$$\H_{2\times2}=\mbox{Her}_2+i\mbox{Her}_2^+\supset0+i\mbox{Her}_2^+
\supset0+i\{z\in\mbox{Her}_2^+\mid\det z=1\}\cong\H^3_\RR,$$
where $\mbox{Her}_2$ stands for the space of $2\times2$--Hermitian
matrices, and $\mbox{Her}_2^+$ that of positive definite ones.
We restrict the groups $\Delta$ and $\Delta(1+i)$ onto the subvariety 
$\H^3_\RR\subset\H_{2\times2}$; we have
$$\Delta|_{\H^3_\RR}=\Gamma,\quad \Delta(1+i)|_{\H^3_\RR}=\Gamma(1+i).$$
Since the Whitehead-link-complement group $W$ is a subgroup of $\Gamma$
of finite index, and so $W$ is commensurable with
$\Gamma(1+i)$,  one can expect that the restrictions of the theta functions 
in the Hermitian setup would give the desired projection $\pi\co  \H^3_\RR\to S^3-L.$
\subsection{Toward an embedding of $\H^3_\RR/W$}\label{sec:4.5}
The expectation stated at the end of the previous section becomes true;
we can construct automorphic functions (out of the theta functions) $f_j$
defined in $\H^3_\RR$ invariant under $W$, and get a map
$$\H^3_\RR\ni x\longmapsto (f_j(x))\in \R^N,\qquad N=13,$$ which gives
an embedding of $\H^3_\RR/W$ onto a semi-algebraic subset of $\R^N$
homeomorphic to $S^3-L$; refer to our paper \cite{MNY} for precise
formulation and proofs. Here we describe what is going on.
\subsubsection{Base camp}
The miraculously lucky fact is that the theta functions
$\theta{a\choose b}(\tau)$ for $a,b\in (\Z[i]/(1+i))^2$ are real
valued on $\H^3_\RR$. (We do not know yet a basic reason for this.)

\begin{figure}[ht!]
\cl{\raise20pt\hbox{\includegraphics[width=4.5cm]{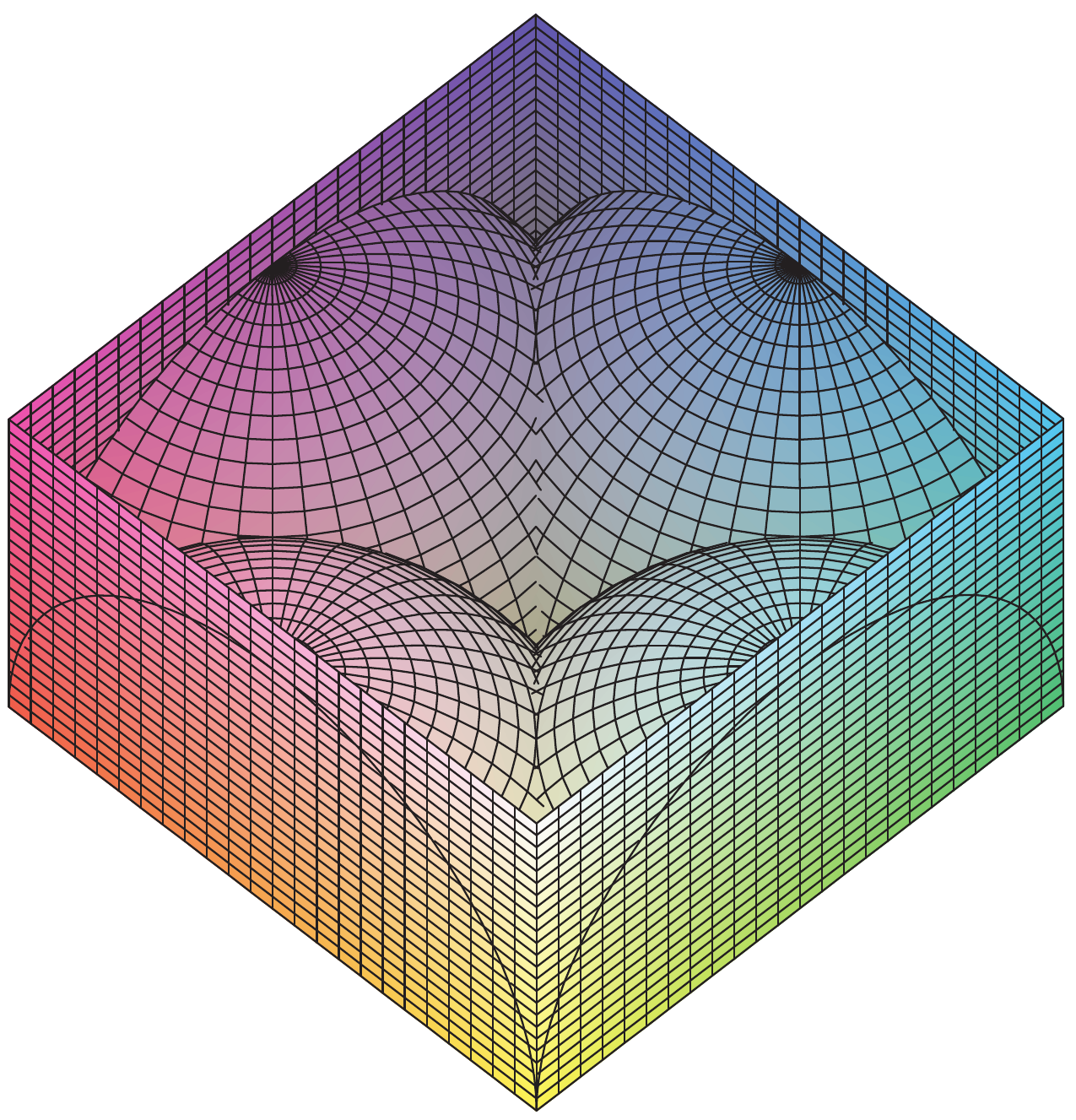}}\qquad
\includegraphics[scale=.9]{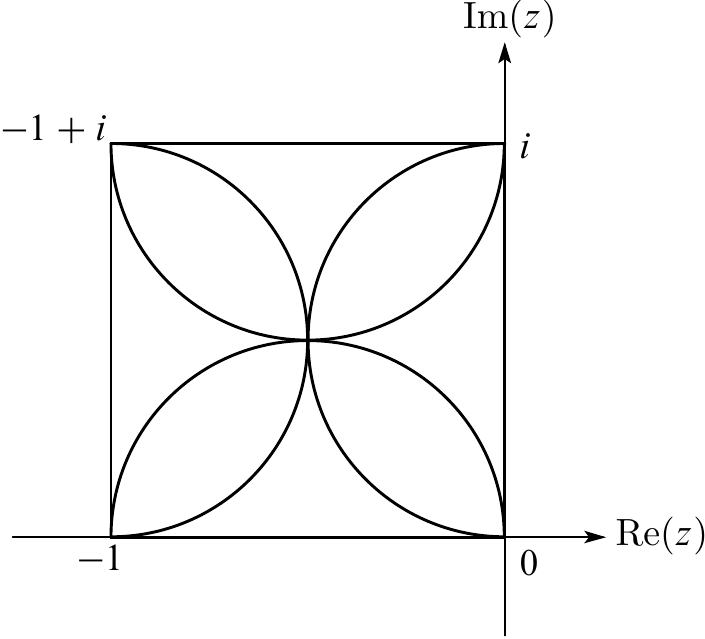}}
\caption{Weyl chamber of $\G^T(2)$}
\label{Weylchamber}.
\end{figure}

Let $\Gamma(2)$ be the principal congruence subgroup of $\Gamma$ of
level 2. The group $\Gamma^T(2):=\langle\Gamma(2),T\rangle$ is a
(hyperbolic) Coxeter group with (ideal) octahedral Weyl chamber (see
\fullref{Weylchamber}). For suitable four such thetas $x_0,\dots,x_3$,
this Weyl chamber ($\sim\H^3_\RR/\Gamma^T(2)$) is mapped
diffeomorphically by the map
$$\H^3_\RR\ni(t,z)\longmapsto\dfrac1{x_0}(x_1,x_2,x_3)\in\R^3,$$
onto the euclidean octahedron
$$\{(t_1,t_2,t_3)\in\R^3\mid |t_1|+|t_2|+|t_3|\le1\}$$
minus the six vertices. Note that this embedding of $\H^3_\RR/\Gamma^T(2)$ is of codimention 0. This embedding plays a role of a base camp of our orbifold-embedding-tour.

\subsubsection{Symmetries of $L$}
Though the codimension of the embedding announced in the beginning of this
subsection (\ref{sec:4.5}) is high ($N-3=10$), we can nevertheless
{\it see} what is happening as follows. We consider a chain of supgroups
of $W$:
$$W\ \overset2\subset\ W_1\ \overset2\subset\  W_2\ \overset2\subset\ \Lambda:=\langle\Gamma^T(2),W\rangle,$$
 each inclusion being of index 2, such that the double coverings
$$(S^3-L\sim)\ \H^3_\RR/W\ \overset{p_1}\longrightarrow\ 
\H^3_\RR/W_1\ \overset{p_2}\longrightarrow\ 
\H^3_\RR/W_2\ \overset{p_3}\longrightarrow\ 
\H^3_\RR/\Lambda$$
of the orbifolds have clear geometric interpretations. 

\begin{figure}[ht!]
\begin{center}
\labellist\small
\pinlabel $F_1$ [r] at 109 45
\pinlabel $F_2$ [t] at 11 226
\pinlabel $F_3$ [t] at 420 226
\pinlabel $L_\infty$ [bl] at 166 380
\pinlabel $L_0$ [bl] at 345 309
\pinlabel 2 [r] at 107 345
\pinlabel $\sq$ [b] at 107 467
\endlabellist 
\includegraphics[width=.5\hsize]{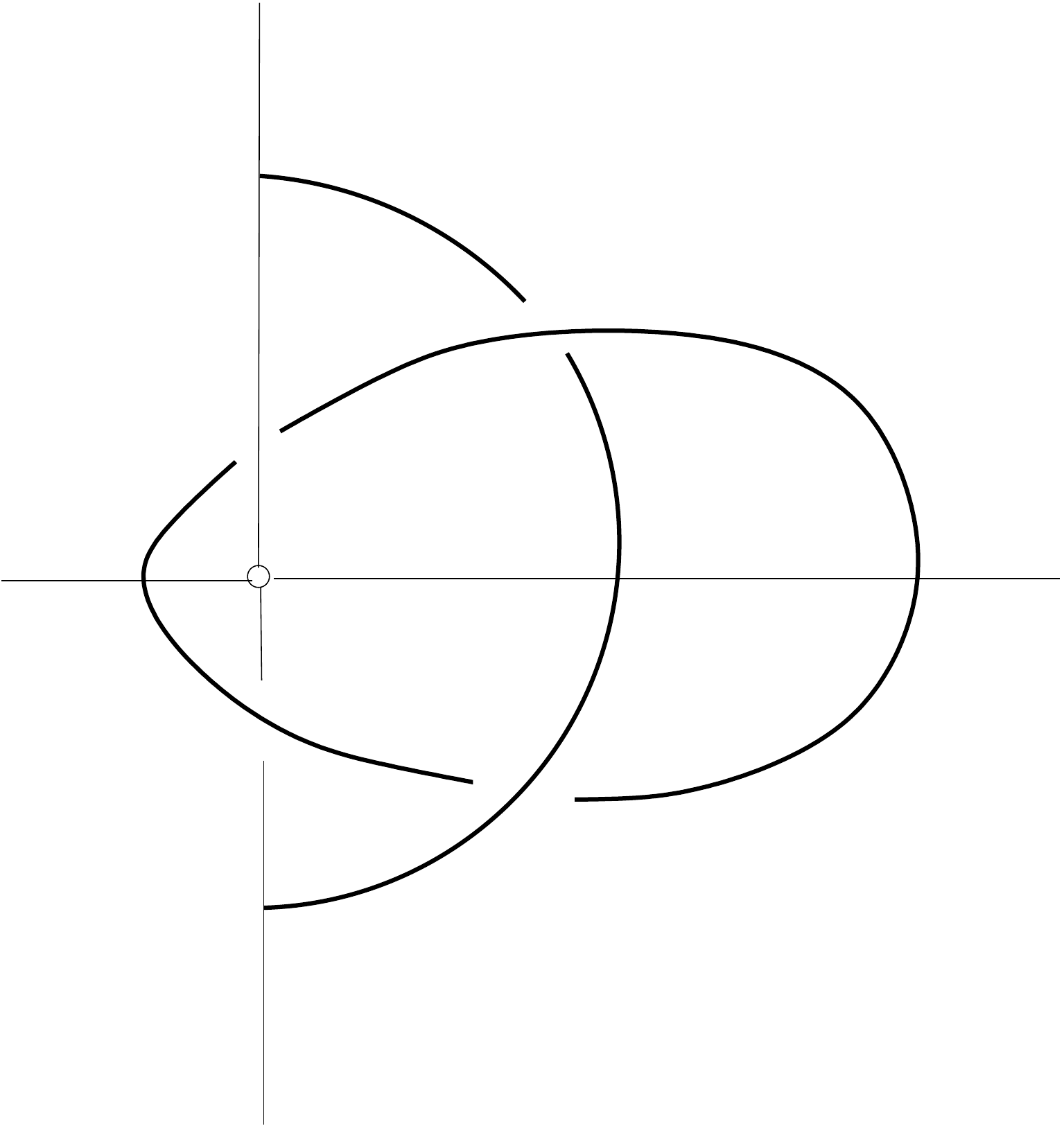}
\end{center}
\caption{The orbifold $\H^3_\RR/W_1$\label{W1}}
\end{figure}

\begin{figure}[ht!]
\begin{center}
\labellist\small
\pinlabel $F_1$ [r] at 107 121
\pinlabel $F_2$ [t] at 23 3
\pinlabel $F_3$ [t] at 323 4
\pinlabel $L_\infty$ [bl] at 167 159
\pinlabel $L_0$ [bl] at 339 93
\pinlabel $\sq$ [b] at 107 245
\pinlabel 2 [r] at 107 210
\pinlabel 2 [r] at -1 3
\pinlabel 2 [b] at 424 5
\endlabellist 
\includegraphics[width=.5\hsize]{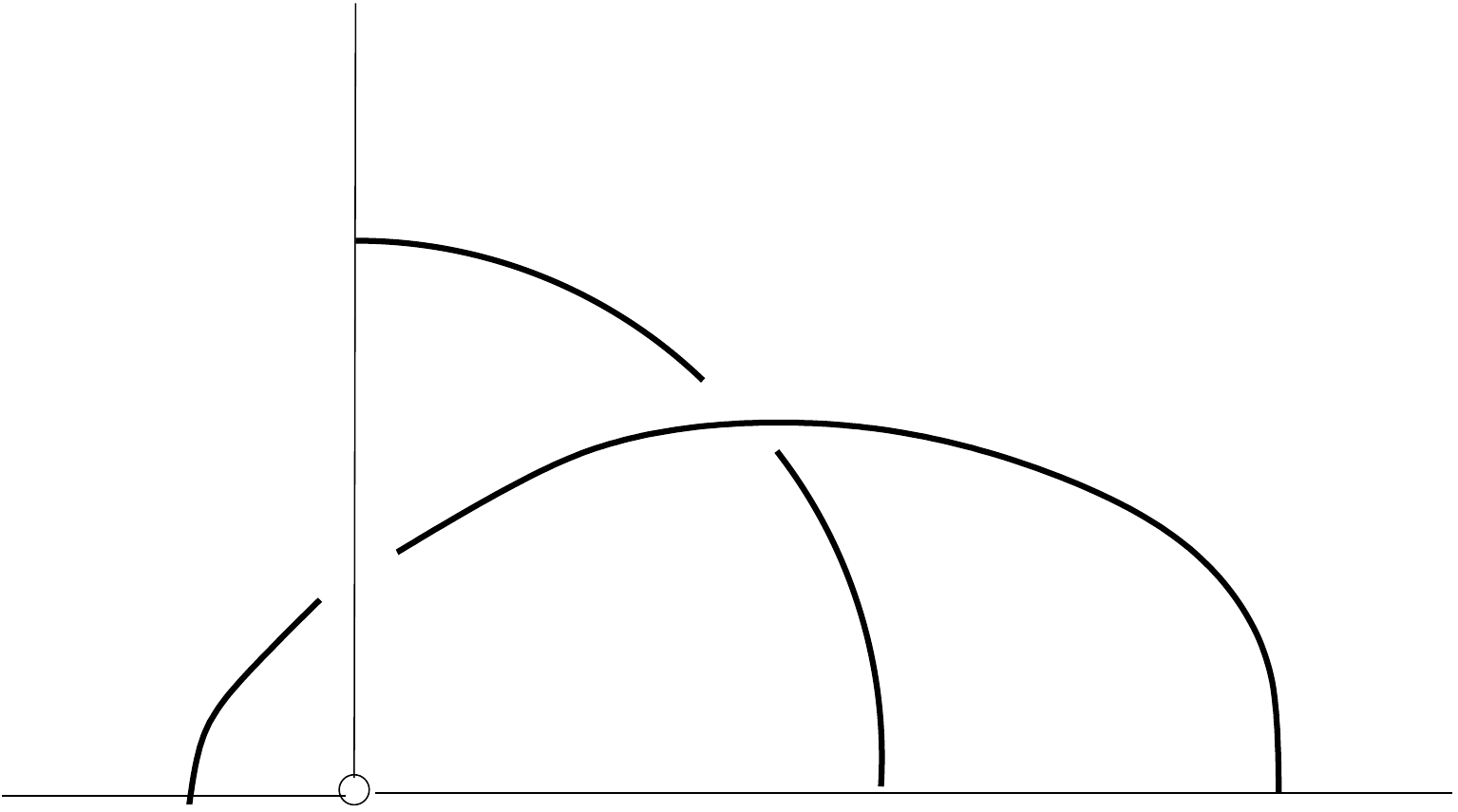}
\end{center}
\caption{The orbifold $\H^3_\RR/W_2$\label{W2}}
\end{figure}

\begin{figure}[ht!]
\begin{center}
\includegraphics[scale=.9]{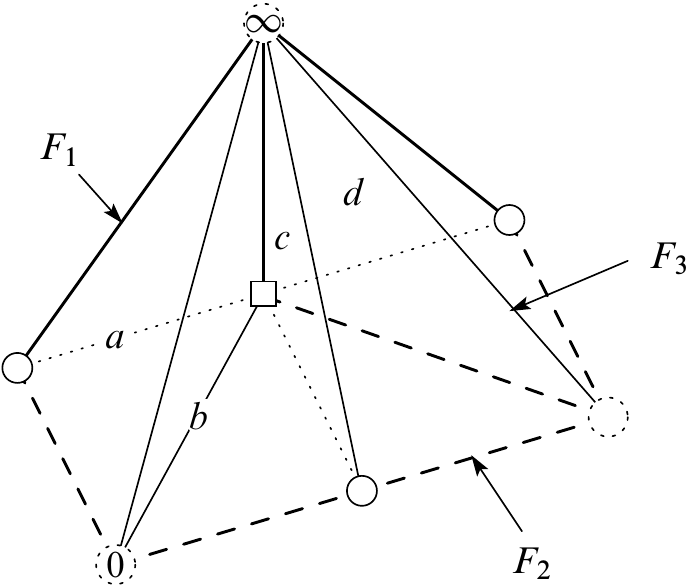}
\end{center}
\caption{A fundamental domain for $W_2$\label{betterpicture}}
\end{figure}

\begin{figure}[ht!]
\begin{center}
\labellist\small
\pinlabel $F_1$ [r] at 50 400
\pinlabel $F_1$ [r] at 126 169
\pinlabel $F_1$ [l] at 529 109
\pinlabel $F_2$ [t] at 4 282
\pinlabel $F_2$ [t] at 394 282
\pinlabel $F_3$ [t] at 267 282
\pinlabel $F_3$ [t] at 611 282
\pinlabel $L_\infty$ [b] at 84 447
\pinlabel $L_\infty$ [tr] at 162 115
\pinlabel $(L_\infty)$ [t] at 525 278
\pinlabel $(L_\infty)$ [l] at 533 143
\pinlabel $L_0$ [bl] at 272 376
\pinlabel $L_0$ [bl] at 638 376
\pinlabel $\sq$ [b] at 126 209
\pinlabel {wall a} at 241 433
\pinlabel {wall b} at 439 433
\pinlabel {wall c} at 50 117
\pinlabel {wall d} at 630 117
\endlabellist 
\includegraphics[width=.9\hsize]{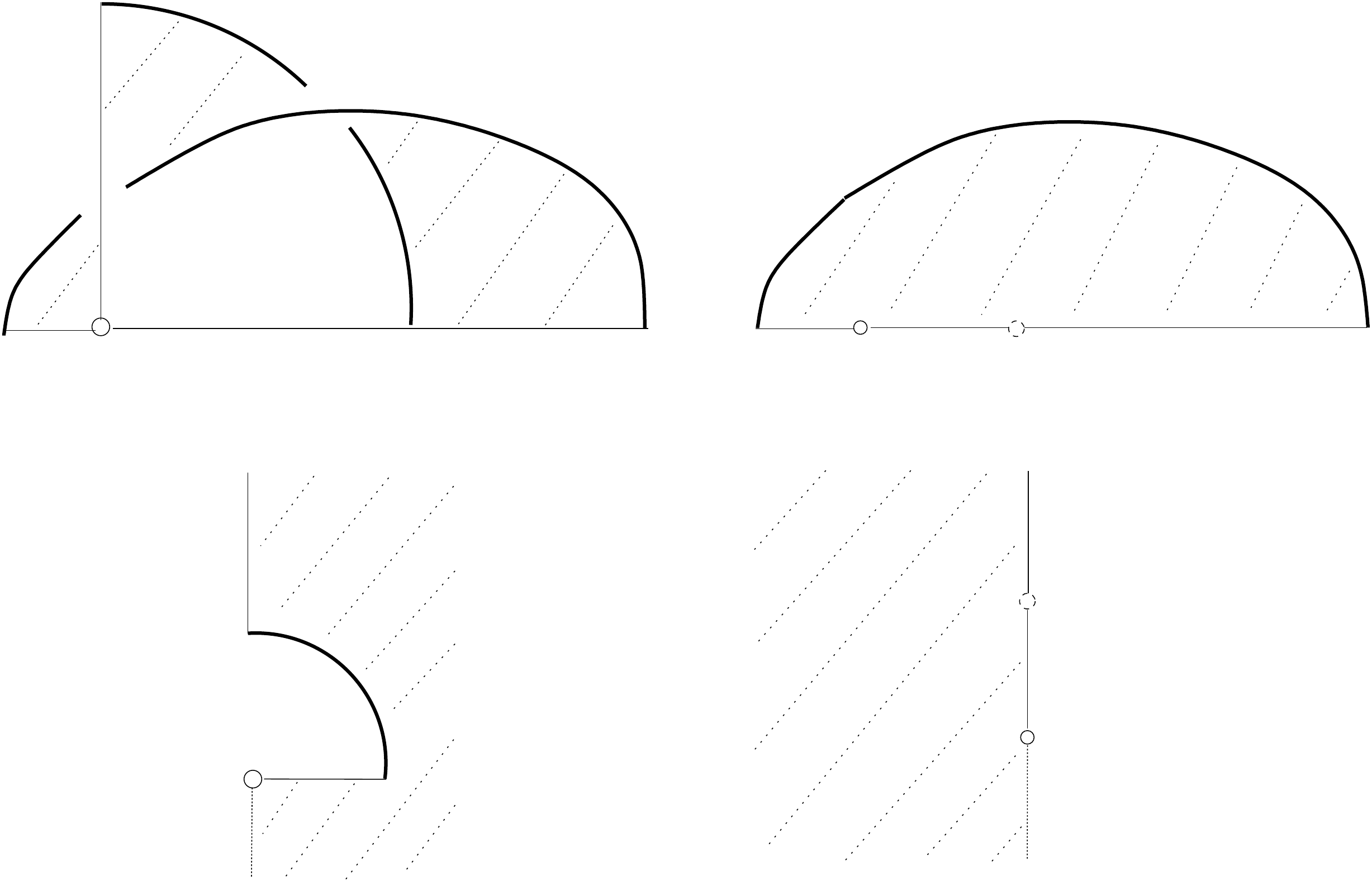}
\end{center}
\caption{The mirror of the reflection in the orbifold
$\H^3/W_2$ is shown as the union of four parts\label{mirror}.}
\end{figure}

\begin{itemize}
\item $p_1$ is the projection modulo the order-2--rotation with axis $F_1\subset  \H^3_\RR/W\sim S^3-L$ (cf Figures \ref{figwhlink}, \ref{W1}),
\item $p_2$ is the projection modulo the order-2--rotation with axis $p_1(F_2)\cup p_1(F_3)\subset \H^3_\RR/W_1$ (cf Figures \ref{W1}, \ref{W2}),
\item $p_3$ is the projection modulo an order-2--reflection with a mirror ($\sim S^2$) in  $\H^3_\RR/W_2$ (cf Figures \ref{W2}, \ref{betterpicture}). The mirror  is shown in \fullref{mirror} as the union of four {\it triangles} labeled as $a,b,c$ and $d$. 
\end{itemize}
The orbifold $\H^3_\RR/\Lambda$ is homeomorphic to a 3--ball with two holes, corresponding to the two strings $L_0$ and $L_\infty$ forming $L$. We have a very simple embedding of $\H^3_\RR/\Lambda$ of codimension 1, as we see in \fullref{sec:4.5.3}.

\subsubsection{Climb up and down}\label{sec:4.5.3}
The group $\Gamma^T(2)$ is a normal subgroup of $\Lambda$, and $\Lambda/\Gamma^T(2)$ is isomorphic to the dihedral group of order eight. To climb up from $\Gamma^T(2)$ to $\Lambda$ is easy. We have only to find rational functions in $x_0,\dots,x_3$ invariant under the dihedral group. Actuallly, by the map
$$\H^3_\RR\ni(z,t)\longmapsto
(\lambda_1,\dots,\lambda_4)=(\xi_1^2+\xi_2^2,\xi_1^2\xi_2^2,\xi_3^2,\xi_1\xi_2\xi_3)\in\R^3,$$
where $\xi_i=x_i/x_0$, the orbifold $\H^3_\RR/\Lambda$ is embedded
into a subdomain of the quadratic hypersurface
$\lambda_2\lambda_3=\lambda_4^2.$ 

To climb down from $\Lambda$ to $W_2$, we need functions which
separate the two sheets of the double cover
$p_3\co \H^3_\RR/W_2\to\H^3_\RR/\Lambda$. These functions are made from
the theta functions $\theta{a\choose b}(\tau)$ with deeper
characteristics $a,b\in (\Z[i]/2)^2$. To determine the image of the
embedding, we make use of the transformation formulas and algebraic
relations among these thetas obtained by Matsumoto \cite{Mat}.

When we climb down from $W_2$ to $W_1$, and finally from $W_1$ to $W$,
we need several functions made by such thetas. Anyway, by functions
made by thetas, we can thus embed $\H^3_\RR/W$ into $\R^{13}$.

\bibliographystyle{gtart}
\bibliography{link}

\end{document}